\documentclass[letterpaper,12pt,leqno]{article}

\usepackage{fullpage}

\usepackage{hyperref}

\usepackage{amsmath}
\usepackage{amsfonts}
\usepackage{amssymb}
\usepackage{tabularx}
\usepackage{color}
\usepackage{fancyhdr}
\usepackage{tikz-cd}
\usepackage{enumitem}

\def\theoremparent{section}
\usepackage{notes}

\numberwithin{equation}{theorem}

\def\A{\mathbb{A}}

\def\cF{\mathcal{F}}

\def\d{\partial}

\def\Der{\mathrm{Der}}

\def\g{\mathfrak{g}}

\def\gr{\mathop{\mathrm{gr}}}

\def\h{\mathfrak{h}}
\def\rmht{\mathrm{ht}}

\def\Ind{\mathop\mathrm{Ind}\nolimits}

\def\k{\mathbf{k}}
\def\kk{\mathfrak{k}}

\def\proof{\noindent{\em Proof:}\ }

\def\qed{\hfill\lower 1em\hbox{$\square$}\vskip 1em}

\def\SS{\mathbb{S}}

\def\sk{\mathop{\mathrm{sk}}}

\def\sVect{\mathrm{sVect}}

\def\toto{\text{\ \raise0.2em\hbox to 0pt {$\to$}\lower0.2em\hbox{$\to$}\ }}

\def\V{\mathcal{V}}

\begin{document}

\title{Gelfand-Fuchs cohomology for affine superspaces $\A^{n,1}$}
\author{Slava Pimenov}
\date{\today}
\titlepage
\maketitle

\tableofcontents
\vskip 5em

\setcounter{section}{-1}
\section{Introduction}
Let $M$ be a differentiable (or algebraic) manifold, and denote by $\V_M = \Gamma(M, TM)$ the space of global section of the tangent bundle.
The Lie bracket of vector fields on $M$ equips $\V_M$ with the structure of Lie algebra, and we are interested in the
cohomology $H^\bullet(\V_M, \k)$. As is often the case with cohomology theories the calculation of $H^\bullet(\V_M, \k)$ can be
split into two steps.
\begin{enumerate}[label=\alph*)]
\item Local calculation of $H^\bullet(\V_{\hat\A^n}, \k)$, where $\hat\A^n$ is a formal $n$-dimensional disk.
\item Some kind of local to global construction.
\end{enumerate}

\begin{nparagraph}
The local calculation is a classical result and can be found in \cite{Fuks}. Consider the topological group $GL(n, \CC)$, and let
$BGL(n)$ be its classifying space. Denote by $p\from EGL(n) \to BGL(n)$ the tautological principal $GL(n)$-bundle over the classifying
space. We will also write $\sk_d BGL(n)$ for the $d$-dimensional skeleton of $BGL(n)$, i.e. the subspace formed by all cells of
dimension up to $d$. Let $X_d = p^{-1}(\sk_d BGL(n))$ be the fiber of the tautological bundle over the $d$-dimensional skeleton.
\end{nparagraph}

\begin{atheorem}[Gelfand-Fuchs]
\label{athm_classic}
We have an isomorphism
$$
H^\bullet(\V_{\hat\A^n}, \k) \ \isom\  H^\bullet(X_{2n}, \k).
$$
\end{atheorem}

The local to global construction can be produced by considering the spectral sequence associated to the diagonal filtration on the
cochain complex $C^\bullet(\V_M, \k)$, which is formed by cochains supported on various diagonal subspaces of $M^q$ and filtered by
inclusion of diagonals (\cite{Gel-Fuchs}). A more modern approach to this construction uses the formalism of factorization algebras
(\cite{Kapranov}).

\begin{nparagraph}[Supermanifolds.]
Now let $M$ be a (smooth or algebraic) supermanifold, in this case even the local picture of the completion of the affine superspace
$\hat\A^{m,n}$ is not fully understood. For convenience
we will write $\V_{m,n} = \V_{\hat\A^{m,n}}$. Some partial results have been obtained regarding the cohomology $H^\bullet(\V_{m,n}, \k)$
(\cite{Fuks-super}) as well as some results concerning cohomology with other coefficients and cohomology of contact vector fields (\cite{BAF}, \cite{FK}).

Let us recall what is known about cohomology of $\V_{m,n}$. The first part of the next theorem was announced without proof in \cite{FL} and
independently proved by Koszul in \cite{Koszul}. In its entirety it was shown in \cite{Fuks-super}.

\end{nparagraph}

\begin{atheorem}[Astashkevich-Fuchs]
\label{athm_mlessn}
For Lie superalgebras $\V_{m,n}$ with $m < n$ we have
\begin{enumerate}[label=\alph*)]
\item $H^\bullet(\V_{0,n}, \k) \isom H^\bullet(S^{2n-1}, \k)$,
\item the inclusion $\V_{0,n} \into \V_{m,n}$ induces an isomorphism $H^\bullet(\V_{m,n}, \k) \isom H^\bullet(\V_{0,n}, \k)$.
\end{enumerate}
\end{atheorem}
Here $S^{2n-1}$ denotes the $(2n-1)$-dimensional sphere. The first part is established using a similar argument to that of the classical case of $\V_{m,0}$,
and the second part is established by comparing cochain complexes for $\V_{0,n}$ and $\V_{m,n}$.
Furthermore, this method can be pushed one step further to obtain the non-trivial result.

\begin{atheorem}[Astashkevich-Fuchs]
\label{athm_nn}
We have an isomorphism
$$
H^\bullet(\V_{n,n}, \k) \ \isom\ H^\bullet(\SS^{2n} GL(n, \CC), \k).
$$
\end{atheorem}
Here $\SS$ denotes the topological suspension functor. At this point it is unclear if this method can be used to continue further
and include cases of $\V_{m,n}$ with $m > n$.
Using a different approach, more in line with the standard calculation used in the classical case, we calculate the cohomology for $\V_{n,1}$.
The main result of this paper is the following theorem.

\begin{atheorem}
\label{athm_n1}
We have an isomorphism
$$
H^\bullet(\V_{n,1}, \k) \ \isom\ H^\bullet(\SS^2 X_{2(n-1)}, \k).
$$
\end{atheorem}

Comparing these results we are led to the following conjectural form for the general case. This conjecture is well known among the specialists
but doesn't seem to appear in the literature.

\begin{aconjecture}
For Lie superalgebras $\V_{m,n}$ with $m, n \ge 0$ we have
$$
H^\bullet(\V_{m,n}, \k) \ \isom\ H^\bullet(\SS^{2n} X_{2(m-n)}, \k).
$$
\end{aconjecture}
Here as before $X_{2(m-n)}$ is the restriction of the tautological principal bundle over $BGL(m)$ to the $2(m-n)$-dimensional skeleton.
Clearly, for $n = 0$ and $n = 1$ the conjecture matches the theorems \ref{athm_classic}, \ref{athm_n1}. For $m = n$ we observe that
the $0$-skeleton $\sk_0 BGL(n)$ is a point and the fiber of the tautological fibration $p$ is $GL(n)$, hence we recover statement
of theorem \ref{athm_nn}. In the case $m < n$ the skeleton $\sk_{2(m-n)} BGL(m)$ is an empty set, and we put by definition
the suspension of an empty set to be the zero-dimensional sphere $S^0$. Therefore, $\SS^{2n} X_{2(m-n)} = S^{2n-1}$ and we
arrive to the statement of theorem \ref{athm_mlessn}.

\vskip 1em
\begin{nparagraph}[Outline of the proof.]
Similar to the classical case we begin the proof by considering Lie subalgebra $\gl(m,n) \into \V_{m,n}$ and constructing the associated
spectral sequence starting from cohomology of $\gl(m,n)$ and converging to the cohomology of $\V_{m,n}$. This is done in the beginning
of section \ref{sec_Vn1}. We arrive to the question of describing the space of invariants in the tensor product
$$
\Lambda^p(V) \tensor \Lambda^p(S^2(V^*) \tensor V).
$$
In the classical case Gelfand and Fuchs introduce multiplicative generators of this space of invariants using graph description. Namely,
they define invariants $\psi_r$ represented by the graphs of the following type.

\vskip 1em
\centerline{
\def\u{0.15em}
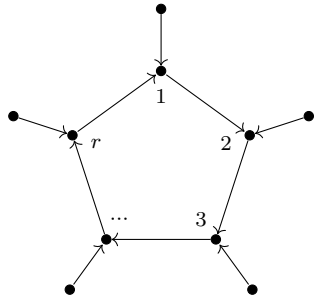
\begin{tikzpicture}
\foreach \i in {0,1,2,3,4}
{
\draw [fill] (90-72*\i:3em) circle [radius=\u] node (p\i) [inner sep=\u] {};
\draw [fill] (90-72*\i:5em) circle [radius=\u] node (q\i) [inner sep=\u] {};
\draw [->] (q\i) -- (p\i);
}
\foreach \i [remember=\i as \j (initially 4)] in {0,1,2,3,4} \draw [->] (p\j) -- (p\i);
\foreach \i in {1,2,3} \node at (90 + 72 - \i*72:2.2em) {$\scriptstyle \i$};
\node at (90 + 72 - 4*72:2.2em) {$\scriptstyle \cdots$};
\node at (90 + 72 - 5*72:2.2em) {$\scriptstyle r$};
\end{tikzpicture}
}
\vskip 1em

Here vertices with only one output correspond to factors $V$ and vertices with two inputs and one output to factors $S^2(V^*) \tensor V$.
Of course, the same $\psi_r$ could be taken as generators of invariants in the case of $\V_{m,n}$ as well, however this description is
ill suited to describe the relations between them. Instead, we adopt a different approach using Schur functors $\Sigma^\lambda$. Each
invariant comes from the one-dimensional invariant subspace $\k \into \Sigma^\lambda(V) \tensor \Sigma^\lambda(V^*)$, for some diagram
$\lambda$ and we identify all diagrams that contribute to the space of invariants.

The next difficulty to resolve is the fact that unlike the classical case, the category of representations of $\gl(m,n)$ in weight $0$
is not semisimple. Therefore, it is not enough just to know the space of invariants in order to calculate the first layer of the spectral
sequence, we need to know exactly what diagram $\lambda$ the invariant came from. Here we restrict our attention to the case $\gl(n,1)$
and perform the necessary calculations in section \ref{sec_gln1}. This is a major difficulty of the super-case compared to the classical one.
In order to apply the methodology of this paper to $\V_{m,n}$ for general $m \ge n$ we need to have an analog of lemma \ref{lemma_gl11}
for $\gl(n,n)$.

Finally, equipped with these two results we complete the proof by analyzing the spectral sequence associated to the Lie subalgebra
$\gl(n,1) \into \V_{n,1}$.

\vskip 1em
The author would like to thank Dmitry Fuchs, Mikhail Kapranov and Vera Serganova who looked at the preliminary version of the paper
and provided valuable remarks.
\end{nparagraph}

\vskip 5em

\section{Recollections}
In this section we introduce notations and recall some basic facts that will be used in this paper. We denote by $\k$ an
algebraically closed field of characteristic $0$.

\begin{nparagraph}[Schur functors.]
Let $S_n$ be a symmetric group on $n$ elements, and $\k[S_n]$ its group algebra. The irreducible representations of $S_n$ are parametrized
by the conjugacy classes of $S_n$, and each such class corresponds to an unordered partition $\lambda = (\lambda_1, \ldots, \lambda_k)$ of $n$,
with $\lambda_1 \ge \lambda_2 \ge \ldots \ge \lambda_k > 0$ and $\sum_{i=1}^k \lambda_i = n$. If we write a representative of the conjugacy class
as a product of cycles, then the integers $\lambda_i$ are the lengths of the cycles. We will denote the corresponding irreducible $S_n$-module
by $M_\lambda$.

It is often convenient to write such a partition $\lambda$ as a Young diagram consisting of $k$ rows such that row $i$ contains $\lambda_i$ cells.
We will call the number $k$ the height of the diagram $\lambda$ and denote it by $\rmht(\lambda)$. Furthermore,
we will write $\lambda'$ for the transposed Young diagram, i.e. a diagram obtained from $\lambda$ by reflecting it along the diagonal. More formally,
we put $\lambda'_j = \max \{ i \mid \lambda_i \ge j\}$.

Consider a symmetric monoidal category $(\C, \tensor)$, i.e. a monoidal category such that for any two objects $V, W \in \C$ we have natural isomorphisms
$c_{V, W}\from V \tensor W \to W \tensor V$ compatible with the monoidal structure and such that the composition $c_{W,V} \circ c_{V,W} = \id_{V \tensor W}$.
Using isomorphisms $c_{V,V}$ the $n$-fold tensor product $V^{\tensor n}$ is equipped with the action of $S_n$ and for each partition $\lambda$
the Schur functor $\Sigma^\lambda$ is defined by
$$
\Sigma^\lambda(V) = V^{\tensor n} \tensor_{\k[S_n]} M_\lambda.
$$

We also denote by $S^n$ and $\Lambda^n$ the symmetric and exterior powers respectively and recall that we have isomorphisms of
$(\gl(V) \times \gl(W))$-modules
$$
S^n(V \tensor W) = \bigoplus_{|\lambda| = n} \Sigma^\lambda(V) \tensor \Sigma^\lambda(W),\quad\quad
\Lambda^n(V \tensor W) = \bigoplus_{|\lambda| = n} \Sigma^\lambda(V) \tensor \Sigma^{\lambda'}(W).
$$

\end{nparagraph}

\begin{nparagraph}[Super vector spaces.]
\label{par_super_vect}
A super vector space over $\k$ is a $\Z/2$-graded vector space $V = (V_0, V_1)$. We denote by $\sVect$ the symmetric monoidal category of
such spaces with the tensor product given by
$$
V \tensor W = (V_0 \tensor W_0 \oplus V_1 \tensor W_1, V_1 \tensor W_0 \oplus V_0 \tensor W_1),
$$
and the symmetry map given by the Koszul sign rule, i.e. for any two homogeneous elements $v \in V$ and $w \in W$ of degrees $\bar v$, $\bar w$ respectively
we have
$$
c_{V, W}(v \tensor w) = (-1)^{\bar v \bar w} w \tensor v.
$$

The dimension of super vector space $V$ is a pair of integers $(m, n)$, where $m = \dim V_0$ and $n = \dim V_1$. The Schur functor $\Sigma^\lambda (V)$
is non-zero if and only if the diagram $\lambda$ is contained in a thick hook with $m$ rows and $n$ columns, in other words $\lambda_i \le n$ whenever
$i > m$.

\end{nparagraph}

\begin{nparagraph}[Lie superalgebra $\gl(m,n)$.]
A Lie superalgebra $\g$ is a Lie algebra in the symmetric monoidal category $\sVect$. Explicitly, it is a super vector space equipped with
the bracket satisfying the commutation relation $[x, y] = -(-1)^{\bar x \bar y}[y, x]$ and the Jacobi relation
$$
[[x, y], z] = [x, [y, z]] - (-1)^{\bar x \bar y} [y, [x, z]],
$$
for any three homogeneous elements $x, y, z \in \g$.

Let $V$ be a super vector space of dimension $(m, n)$, we denote by $\gl(m,n)$ the Lie superalgebra of all $\k$-endomorphisms of $V$
(that do not necessarily preserve the $\Z/2$ grading). Let us fix a basis $\{v_i\}$ of $V$ and write $e_{ij} = v_i \tensor v_j^*$
for the corresponding basis of $\gl(m,n)$. Clearly, we have $\bar e_{ij} = \bar v_i + \bar v_j \in \Z/2$.

This Lie superalgebra is equipped with the non-degenerate invariant super-symmetric bilinear form, defined by
$$
(e_{ij}, e_{ji}) = (-1)^{\bar v_i},
$$
and with all other pairings equal to $0$. Denote by $\h \subset \gl(m,n)$ the abelian subalgebra of diagonal matrices, the restriction
of this bilinear form to $\h$ is also non-degenerate, and we use it to identify $\h$ with $\h^*$. Let $\alpha \in \h$ be a root of
$\gl(m,n)$, we denote $H_\alpha$ the hyperplane in $\h$ orthogonal to $\alpha$. We say that the root $\alpha$ is isotropic if $\alpha \in H_\alpha$.

Consider the group $S_m \times S_n$ acting on the space of functions $\k[\h]$ by permuting the first $m$ and the last $n$ coordinates.
Then we have the following description of the center $Z\gl(m,n)$ of the universal enveloping algebra $U\gl(m,n)$ using the Harish-Chandra map $\psi$.

\end{nparagraph}

\begin{theorem}[\cite{Musson}]
The map $\psi\from Z\gl(m,n) \to \k[\h]$ is a monomorphism and the image consists of $(S_m \times S_n)$-invariant functions $f$ such that
$
f(\lambda) = f(\lambda + t \alpha)
$
for any $t \in \k$, whenever $\lambda \in H_\alpha$ for an isotropic root $\alpha$.
\end{theorem}

\begin{nparagraph}[Cohomology of $\gl(m,n)$.]
\label{par.cohom.glmn}
Let $\g = \gl(m,n)$, and $V$ as before the standard representation of $\g$. We will be interested in cohomology spaces
$H^\bullet(\g, \Sigma^\alpha V \tensor \Sigma^\beta V^*)$. First, denote by $h$ the identity element of $\g$. Since the Lie algebra
acts trivially on its cohomology we immediately see that this cohomology space can only be nontrivial if $|\alpha| = |\beta|$.
Furthermore, the Hom-space $\R\Hom_\g(\Sigma^\beta V, \Sigma^\alpha V)$ can be non-zero only if the two $\g$-modules $\Sigma^\alpha V$
and $\Sigma^\beta V$ have the same generalized central character (i.e. supported at the same point in $\Spec Z\g$). However, due to the
description of the center above and the fact that all $\alpha_i$ and $\beta_i$ are non-negative this can only happen if $\alpha = \beta$.

Consider a Lie subalgebra $\g' \subset \g$, and denote by $\Ind_{\g'}^\g \from \g'\Mod \to \g\Mod$ the induction functor
$$
\Ind_{\g'}^\g (M) = \Hom_{U\g'}(U\g, M).
$$

\end{nparagraph}

\begin{theorem}
\label{thm.cohom.glmn}
Let $\g = \gl(m,n)$ and assume $m \ge n$, then
\begin{enumerate}[label=\alph*)]
\item The inclusion $\gl(m) \into \g$ induces isomorphism of cohomology
$$
H^\bullet(\g, \k) \isom H^\bullet(\gl(m), \k) \isom \Lambda[e_1, e_3, \ldots, e_{2m-1}],
$$
where each generator $e_i$ is of cohomological degree $i$.
\item Let $\g_0  = \gl(m) \oplus \gl(n) \into \g$ be the even Lie subalgebra of $\g$, then
$$
H^\bullet(\g, \Ind_{\g_0}^\g (\k)) \isom H^\bullet(\g_0, \k) \isom \Lambda[e_1, e_3, \ldots, e_{2m-1}] \tensor \Lambda[e'_1, e'_3, \ldots, e'_{2n-1}].
$$
\end{enumerate}
\end{theorem}

\proof
The first part is the well known calculation and we refer to \cite{Fuks} for details. The second follows immediately from the definition
of the induction functor and the fact that $U\g$ is a free $U\g_0$-module.
\qed

\begin{nparagraph}[Lie superalgebra $\V_{m,n}$.]
Consider affine superspace $\A^{m,n} = \Spec \O_{\A^{m,n}}$, where
$$
\O_{\A^{m,n}} = \k[x_1, \ldots x_m, \xi_1, \ldots \xi_n]
$$
is the polynomial super-algebra of functions with even coordinates $x_i$ and odd coordinates $\xi_j$. We denote by $\hat\O_{\A^{m,n}}$ its completion
at point $0 \in \A^{m,n}$, i.e. $\hat\O_{\A^{m,n}} = \k[[x_1, \ldots, x_m, \xi_1, \ldots, \xi_n]]$, equipped with the power series topology.
We will be interested in the Lie superalgebra of continuous self derivations
$$
\V_{m,n} = \Der_{\mathrm{cont}} (\hat\O_{\A^{m,n}}).
$$

Explicitly, the Lie superalgebra $\V_{m,n}$ consists of elements of the form $\sum f_i {\d_{x_i}} + \sum g_j {\d_{\xi_j}}$, with the
obvious bracket. It contains the subalgebra $\gl(m,n)$ spanned by linear derivations $\{x_i \d_{x_j}, x_i \d_{\xi_j}, \xi_i \d_{x_j}, \xi_i \d_{\xi_j}\}$.

\end{nparagraph}

\begin{nparagraph}[Spectral sequence.]
Let $\g$ be a Lie superalgebra, $\h \into \g$ a Lie subalgebra and $M$ a $\g$-modules. The chain complex $\Lambda^\bullet \g \tensor M$ is equipped
with an increasing filtrations by the number of elements from $\h$, which induces a decreasing filtration on the cochain complex $C^\bullet(\g, M)$.
We have the spectral sequence associated to this filtration
$$
E_1^{pq} = H^q(\h, \ \Hom(\Lambda^p (\g / \h), M)\ ) \Rightarrow H^{p+q}(\g, M).
$$

This will be the main computational tool in this paper, both for the cohomology of $\gl(n,1)$ as well as for the cohomology of $\V_{n,1}$.

\end{nparagraph}

\vskip 5em
\section{Cohomology of $\gl(n,1)$}
\label{sec_gln1}
Throughout this section we consider $\g = \gl(n,1)$, and $V$ will be the standard representation of $\gl(n,1)$. The section is dedicated to calculating
the cohomology $H^\bullet(\g, \Sigma^\alpha V \tensor \Sigma^\alpha V^*)$. We will proceed by induction on $n$, with the case $n = 1$ being the
base of induction. For $\gl(1, 1)$ the calculation can be done directly, since we have a complete description of the category of $\g$-modules
in this case (\cite{FKV}).

\begin{nparagraph}[Representations of $\gl(1,1)$.]
We begin be recalling relevant details concerning the category of $\gl(1,1)$-modules. As before we put $h \in \g = \gl(1,1)$ to be the identity
element, and put
$$
e = \begin{pmatrix}
0 & 1 \\
0 & 0
\end{pmatrix},
\quad
f = \begin{pmatrix}
0 & 0 \\
1 & 0
\end{pmatrix},
\quad
g = \begin{pmatrix}
1 & 0 \\
0 & -1
\end{pmatrix}.
$$

The category of $\g$-modules in non-zero $h$-weight $\lambda$ is semisimple, and the simple modules are two-dimensional with basis vectors
of $(h, g)$-weights $(\lambda, \gamma)$ and $(\lambda, \gamma - 2)$ and the action of $e$ and $f$ are as follows.

$$
{
\def\bbt#1#2{\mathop{\bullet}\limits_{#1}^{#2}}
\begin{tikzcd}[column sep=5em]
\bbt{\gamma-2}{\lambda} \ar[r, shift left, "e"] & \bbt{\gamma}{\lambda} \ar[l, shift left, "f"]
\end{tikzcd}
}
$$
We write $V_{\lambda,\gamma}^{\wbar 0}$ for the simple module as above if the vector of $g$-weight $\gamma$ is even and $V_{\lambda,\gamma}^{\wbar 1}$
if it is odd. In particular we have $V = V_{1, 1}^{\wbar 0}$.
Consider a hook diagram $\lambda = (c, 1, \ldots, 1)$, of height $r$ (so that $|\lambda| = c + r - 1$). We have
$$
\Sigma^\lambda(V) = V_{c + r - 1, c - r + 1}^{\wbar{r - 1}}.
$$

The tensor product $\Sigma^\lambda(V) \tensor \Sigma^\lambda(V^*)$ is isomorphic to the induced module $\Ind_{\g_0}^{\g} (\k)$.
$$
\begin{tikzcd}
& \bullet \ar[dl, "e"'] \ar[dr, "f"] & \\
\bullet \ar[dr, "f"'] && \bullet \ar[dl, "-e"] \\
& \bullet &
\end{tikzcd}
$$

\end{nparagraph}

\begin{lemma}
\label{lemma_gl11}
Let $\g = \gl(1, 1)$ and $V$ the standard representation of $\g$, then
$$
H^\bullet(\g, \Sigma^\lambda(V) \tensor \Sigma^\lambda(V^*)) = \begin{cases}
\k[e_1],&\text{if $|\lambda| = 0$,}\\
\k[e_1, e'_1],&\text{otherwise},
\end{cases}
$$
where both $e_1$ and $e'_1$ are of cohomological degree $1$.
\end{lemma}

\proof
The statement follows immediately from the discussion above and theorem \ref{thm.cohom.glmn}.

\qed

\begin{proposition}
\label{prop_cohom_gln1}
Let $\g = \gl(n,1)$, and $V$ the standard representation of $\g$, then
$$
H^\bullet(\g, \Sigma^\lambda(V) \tensor \Sigma^\lambda(V^*)) = \begin{cases}
\k[e_1, \ldots, e_{2n-1}],&\text{if $\rmht(\lambda) \le n - 1$,}\\
\k[e_1, \ldots, e_{2n-1}, e'_1],&\text{otherwise},
\end{cases}
$$
where generators $e_i$ are of cohomological degree $i$ and $e'_1$ is of degree $1$.
\end{proposition}

\proof
The case $n = 1$ is the statement of the previous lemma. Assume that $n > 1$ and that the proposition has been established
for all $\gl(m,1)$ with $m < n$. Consider the Lie subalgebra $\h = \gl(1) \oplus \gl(n - 1, 1) \into \g$ embedded as diagonal blocks.
We will write $E$ for the standard representation of $\gl(1)$ (it is of course a one-dimensional even vector space, however we use letter $E$ instead of
$\k$ to distinguish it from the trivial representation), and $W$ for the standard representation of $\gl(n-1,1)$. We have
$$
\g / \h \ \isom\  E \tensor W^* \oplus E^* \tensor W,
$$
and therefore

$$
\Lambda^p(\g/\h) \ \isom\  \bigoplus_{i + j = p} S^i(E) \tensor \Lambda^i(W^*) \tensor S^j(E^*) \tensor \Lambda^j(W).
$$

Furthermore, we have $V = E \oplus W$ and hence
\begin{equation}
\label{equ_VEW}
\Sigma^\lambda(V) \ \isom\  \bigoplus_\beta S^k(E) \tensor \Sigma^\beta(W),
\end{equation}
where the sum is taken over all subdiagrams $\beta \subset \lambda$, such that the skew diagram $\lambda \backslash \beta$
can be collapsed into a single row diagram, and $k$ is the size of $\lambda \backslash \beta$.

Consider the spectral sequence for the subalgebra $\h \into \g$:
$$
E_1^{pq} = H^q(\h, \Hom(\Lambda^p(\g / \h), \Sigma^\lambda(V) \tensor \Sigma^\lambda(V^*))) \Rightarrow
H^{p+q}(\g, \Sigma^\lambda(V) \tensor \Sigma^\lambda(V^*)).
$$
Using the decompositions described above we see that the first layer of this spectral sequence consists of terms of the form
\begin{equation}
\label{equ_terms}
H^\bullet(\gl(1) \oplus \gl(n-1,1), S^k(E) \tensor S^l(E^*) \tensor \Sigma^\beta(W) \tensor \Sigma^\gamma(W^*)).
\end{equation}
From the discussion in paragraph \ref{par.cohom.glmn} it follows that these terms are non-zero only if $\gamma = \beta$ (in which
case we necessarily have $k = l$). We will refer to the sum of all such terms as the component of type $\beta$. For convenience
we introduce notation $\Delta^\lambda(V) = \Sigma^\lambda(V) \tensor \Sigma^\lambda(V^*)$.

\begin{nparagraph}[Case of $\rmht(\lambda) \ge n$.]
For any Young diagram $\alpha$ we can construct diagram $\alpha^+$ by adding one box in the first column. Recall that for the
Schur functor $\Sigma^\lambda(V)$ to be non-zero the diagram $\lambda$ must have no more than $n$ boxes in the second column,
and similarly diagrams $\beta$ in (\ref{equ_VEW}) and (\ref{equ_terms}) no more than $(n-1)$ boxes. Using this we observe that
the diagram combinatorics in the first layer of spectral sequences for modules $\Delta^\lambda(V)$ and $\Delta^{\lambda^+}(V)$
are identical, in the sense that we have a bijection between terms of type $\beta$ in the former and terms of type $\beta^+$ in the latter.
Moreover, since $\rmht(\beta) \ge (n-1)$ and by inductive assumption we have isomorphism of cohomology
$$
H^\bullet(\gl(n-1,1), \Delta^\beta(W)) \isom H^\bullet(\gl(n-1,1), \Delta^{\beta^+}(W))
$$
we see that the two spectral sequences are isomorphic and converge to the same result. Therefore it is enough to consider the
case $\rmht(\lambda) = n$.

We split the first layer $E_1$ into two parts, the first one $E'_1$ will consist of all components of type $\beta$ with
$\rmht(\beta) = n-1$, and the second one $E''_1$ will consist of all the remaining components (i.e. components of type $\beta$ with
$\rmht(\beta) \ge n$).

Consider also Lie subalgebra $\kk = \gl(1) \oplus \gl(n-1) \into \gl(n)$, and let $V_0$ and $W_0$ be the standard representations
of $\gl(n)$ and $\gl(n-1)$ respectively. We have the corresponding spectral sequence
$$
F_1^{pq} = H^q(\kk, \Hom(\Lambda^p(\gl(n) / \kk), \Delta^\lambda(V_0))) \Rightarrow H^{p+q}(\gl(n), \Delta^\lambda(V_0)).
$$
The latter cohomology is isomorphic to the cohomology of the trivial $\gl(n)$-module, and as was mentioned in theorem \ref{thm.cohom.glmn}
it is isomorphic to the exterior algebra $\k[e_1, \ldots, e_{2n-1}]$.

Observe that the diagram combinatorics of the part $E'_1$ and $F_1$ are identical, and moreover by inductive assumption we have
$$
H^\bullet(\gl(n-1,1), \Delta^\beta(W)) \isom H^\bullet(\gl(n-1), \Delta^\beta(W_0)) \tensor \k[e'_1].
$$
So if we show that the second part $E''_1$ is acyclic we will conclude that the spectral sequence $E_\bullet$ converges to
$\k[e_1, \ldots, e_{2n-1}] \tensor \k[e'_1]$ as required.

In fact we will show that every component of type $\beta$ in $E''_1$ is acyclic. The process described in the beginning or the proof
by which diagram $\beta$ is obtained from $\lambda$ consists of choosing a common subdiagram $\alpha$ such that the difference $\lambda \backslash \alpha$
is collapsible to a single row diagram and the difference $\beta \backslash \alpha$ collapses to a single column diagram. We denote
$k = |\lambda \backslash \alpha|$. Of course there could be multiple such ways to obtain $\beta$ from $\lambda$, so let us pick $\alpha$ with minimal $k$
(such $\alpha$ is unique). We call a box in $\alpha$ {\em flippable} if by removing it we obtain a diagram $\alpha_1$ which provides another
way to construct $\beta$ out of $\lambda$. It is easy to see that a box is flippable if and only if it has no boxes to the right in $\beta$
and has no boxes on top of it in $\lambda$. We denote by $d = d(\beta, \lambda)$ the number of flippable boxes. Since $\rmht(\lambda) = n$ and $\rmht(\beta) \ge n$
we see that $d > 0$.

By construction of the spectral sequence the differential in $E_1$ corresponds to the process of flipping boxes, therefore we find that
the component of $E''_1$ of type $\beta$ is of the form
$$
\Delta^\beta(W) \tensor \bigotimes_1^{2d} \left( \k \to \k \right).
$$
The $2d$ factors come from $d$ flippable boxes involving terms $\Sigma^\beta(W)$ and another $d$ for terms $\Sigma^\beta(W^*)$. Since
each factor is acyclic and $d > 0$ we conclude that the component of $E''_1$ of type $\beta$ is acyclic.

This concludes the proof of the second part of the proposition.

\end{nparagraph}

\begin{nparagraph}[Case of $\rmht(\lambda) \le n-1$.]
\label{prop_proof_n-1}
As before we split the first layer of the spectral sequence into two parts $E'_1$ consisting of components of type $\beta$ with $\rmht(\beta) \le n-1$,
and $E''_1$ consisting of components with $\rmht(\beta) \ge n$.

First, we compare $E'_1$ with $F_1$. Using the inductive assumption we see that components of type $\beta$ with $\rmht(\beta) < n - 1$ provide
isomorphic terms in $E'_1$ and $F_1$, while those with $\rmht(\beta) = n - 1$ have extra summand of the form $\k[e_1, \ldots, e_{2n-3}]e'_1 z_\beta$,
where $z_\beta$ is just a basis vector introduced for bookkeeping purposes. Moreover, from the discussion at the end of the last paragraph we
see that these extra terms survive into $E_2$ if and only if $\beta$ has no flippable boxes, i.e. $d(\beta, \lambda) = 0$. We will show that $E''$
together with these extra terms from $E'$ degenerate after $E_2$. Therefore, starting from the third layer spectral sequences $E_\bullet$ and $F_\bullet$
are isomorphic and converge to the same thing $\k[e_1, \ldots, e_{2n-1}]$ as required.

As was discussed before, the only contribution to $E''_2$ comes from diagrams $\gamma$ without flippable boxes and it is of the form
$\k[e_1, \ldots, e_{2n-3}, e'_1] z_\gamma$. It is easy to see that such diagram $\gamma$ can be obtained from some diagram $\beta$ with $\rmht(\beta) = n - 1$
and $d(\beta, \lambda) = 0$ by repeated application of the $+$-operation. For every diagram $\lambda$ we put
$$
l(\lambda) = l_{n-1}(\lambda) = \sum_{i = 1}^{n-1} \lambda_i,
$$
and consider decreasing filtration of $E_2$ by $\cF_{\ge l} E_2$ consisting of all components of type $\beta$ with $l(\beta) \ge l$. In order
to establish the acyclicity of the subcomplex in question it is enough to show acyclicity of the associated graded for this filtration.
Next, consider another grading on $E_2$ by the number of $e'_1$ and the induced grading on the associated graded space $\gr_{\cF} E_2$.
This turns $\gr_{\cF} E_2$ into a collection of bicomplexes where one differential acts on $e_1$ but preserves $e'_1$, while the second differential acts
trivially on generators $e_i$ and $d(e'_1 z_\gamma) = z_{\gamma^+}$. Now, by taking cohomology with respect to the latter differential we see
that the associated graded of $E''_2$ together with the extra terms from $E'_2$ mentioned before is acyclic, and hence the $E''_2$ itself with
those extra terms is also acyclic.

This completes the proof of the proposition.

\qed

\end{nparagraph}

\begin{remark}
This proposition combined with the theorem \ref{thm.cohom.glmn} suggests that the only contribution to the cohomology
$H^\bullet(\g, \Sigma^\lambda(V) \tensor \Sigma^\lambda(V^*))$ comes from the indecomposable modules $\k$ and $\Ind_{\g_0}^{\g}(\k)$.
Indeed, suppose we know that the invariant subspace $\k \into \Sigma^\lambda(V) \tensor \Sigma^\lambda(V^*)$ is either a direct
summand itself or a part of an indecomposable direct summand isomorphic to $\Ind_{\g_0}^{\g}(\k)$. These two possibilities can be
easily distinguished. First, observe that the simple odd root of $\g$ acts acyclically on $\Ind_{\g_0}^{\g}(\k)$. Next, write
$V = X \oplus Y$, where $X$ is the standard representation of $\gl(n-1)$, and $Y$ is the standard representation of $\gl(1, 1)$.
Since $\Sigma^\beta(X) = 0$ whenever $\rmht(\beta) > n - 1$ we see that in the decomposition of $\Sigma^\lambda(V)$ with $\rmht(\lambda) > n - 1$
every summand involves at least one $Y$, and similarly for $\Sigma^\lambda(V^*)$. Therefore, from the representation theory of $\gl(1, 1)$-modules
recalled above it follows that the simple odd root of $\gl(n,1)$ acts acyclically on $\Sigma^\lambda(V) \tensor \Sigma^\lambda(V^*)$, which in turn
implies that it can not contain the trivial representation as a direct summand. Then comparing theorem \ref{thm.cohom.glmn} and the proposition
we conclude that $\Sigma^\lambda(V) \tensor \Sigma^\lambda(V^*)$ with $\rmht(\lambda) \le n - 1$ has $\k$ as a direct summand, and no summands
other than $\k$ and $\Ind_{\g_0}^{\g}(\k)$ can contribute to the cohomology.

It would be interesting to give an independent description of the structure of the module $\Sigma^\lambda(V) \tensor \Sigma^\lambda(V^*)$, as it
would provide alternative (and possibly more meaningful) proof of the proposition \ref{prop_cohom_gln1}.

\end{remark}

\vskip 5em
\section{Cohomology of $\V_{n,1}$}
\label{sec_Vn1}
We begin this section with the standard opening argument used in calculation of $H^\bullet(\V_{m,0}, \k)$ (for details we refer to \cite{Fuks}),
which is also applicable in the general case of $\V_{m,n}$.

Consider Lie subalgebra $\gl(m,n) \into \V_{m,n}$, and let $V$ be the standard representation of $\gl(m,n)$. We have the associated spectral sequence
$$
E_1^{pq} = H^q(\gl(m,n), \Hom(\Lambda^p(\V_{m,n} / \gl(m,n)), \k)) \Rightarrow H^{p+q}(\V_{m,n}, \k).
$$
As a $\gl(m,n)$-module the associated graded $\gr \V_{m,n}$ (recall that $\V_{m,n}$ is the completed Lie superalgebra of vectorfields) is isomorphic to
$$
\gr \V_{m,n} \ \isom\ \bigoplus_{n \ge 0} \left( S^n(V) \tensor V^* \right).
$$
In this decomposition the Lie subalgebra $\gl(m,n) \into \V_{m,n}$ corresponds to the summand $V \tensor V^*$. Therefore, as a $\gl(m,n)$-module
the space $\Hom(\Lambda^p(\V_{m,n} / \gl(m,n)), \k)$ is isomorphic to
\begin{equation}
\label{equ_invar}
\bigoplus_{\sum p_i = p} \Lambda^{p_i} \left(S^i(V^*) \tensor V\right),
\end{equation}
where $i \ge 0$ and $i \neq 1$. In order for the cohomology to be non-zero we need equal number of instances of $V$ and $V^*$ in the above
expression, which leads to
$$
\sum_{i \ge 0, i \neq 1} p_i = \sum_{i \ge 0, i \neq 1} i p_i,
$$
or equivalently, $p_0 = \sum\limits_{i \ge 2} (i-1)p_i$.

On the other hand, as we discussed in \ref{par.cohom.glmn} the contribution to cohomology only comes from terms $\Sigma^\lambda(V) \tensor \Sigma^\lambda(V^*)$,
and the latter contains one-dimensional invariant subspace. Due to the symmetry of $S^i(V^*)$ in (\ref{equ_invar}) in order to obtain an invariant
we must have
$$
p_0 \le \sum_{i \ge 2} p_i.
$$
Combining the two conditions we find that $p_i = 0$ for $i \ge 3$ and $p_0 = p_2$. In other words, it is enough to focus our attention on the
space $\Lambda^p(V) \tensor \Lambda^p(S^2(V^*) \tensor V)$.

\begin{lemma}
\label{lemma_invar}
Let $\g = \gl(m,n)$, and $V$ the standard representation of $\g$, then
$$
H^\bullet(\g, \Lambda^p(V) \tensor \Lambda^p(S^2(V^*) \tensor V)) \ \isom\ 
\bigoplus_{|\lambda| = p} H^\bullet(\g, \Sigma^{\wtilde\lambda}(V) \tensor \Sigma^{\wtilde\lambda}(V^*)),
$$
where $\wtilde\lambda$ is obtained from $\lambda$ by shifting it to the right and adding first column of height $p$. In other words,
$\wtilde\lambda$ is a diagram of height $p$, such that $\wtilde\lambda_i = \lambda_i + 1$ for $1 \le i \le p$.
\end{lemma}

\proof
First, recall that we have
$$
\Lambda^p(S^2(V^*) \tensor V) \ \isom\ \bigoplus_{|\lambda| = p} \Sigma^\lambda(S^2(V^*)) \tensor \Sigma^{\lambda'}(V).
$$
Since the contribution to the cohomology comes only from terms of the form $\Sigma^\alpha(V) \tensor \Sigma^\alpha(V^*)$, we can
disregard everything else. Observe that the product $\Lambda^p(V) \tensor \Sigma^{\lambda'}(V)$ is the sum of $\Sigma^\beta(V)$ with
$\rmht(\beta) \ge p$. On the other hand, the height of diagrams appearing in $\Sigma^\lambda(S^2(V^*))$ is less or equal than $p$.
Hence, it is enough to consider only diagrams of height $p$.

Notice, that in the product $\Lambda^p(V) \tensor \Sigma^{\lambda'}(V)$
there is only one such diagram, namely $\wtilde{\lambda'}$. Furthermore, in constructing diagrams in $\Sigma^\lambda(S^2(V^*))$ the first
box of every pair must go into the first column, and since they anticommute with each other the commutation relation for the rest of the
boxes are opposite to that of $\lambda$, i.e. they form the transposed diagram $\lambda'$. Thus, we have a term $\Sigma^{\wtilde{\lambda'}}(V^*)$
in the expansion of $\Sigma^\lambda(S^2(V^*))$. This completes the proof of the lemma.

\qed

\begin{remark}
\label{rem_invar}
From the lemma it is immediate to see what invariants can appear in
$$
\Lambda^\bullet(V) \tensor \Lambda^\bullet(S^2(V^*) \tensor V)
$$
and give the relations for the multiplicative generators. Specifically, the diagram $\wtilde\lambda$ must be contained in a thick hook
with $m$ rows and $n$ columns (see \ref{par_super_vect}).

As an example, we illustrate it in the cases $\V_{m,0}$ and $\V_{0,n}$, that have been previous handled using a different argument
in \cite{Fuks} and \cite{Fuks-super}. In the case $\V_{m,0}$ we have $\rmht(\wtilde\lambda) \le m$, therefore the space of invariants
is generated by Chern classes $c_i$ of degree $2i$, for $1 \le i \le m$ (corresponding to columns of height $i$), with the condition
that the total degree is less or equal than $2m$ (the height of the first column).

For $\V_{0,n}$ the diagram $\wtilde\lambda$ must be contained in the first $n$ columns, so the space of invariants is generated by
the Segre classes $s_i$ of degree $2i$, for $1 \le i \le n - 1$ with no relations. Here we also notice that since $s_i \equiv c_i$ modulo
subalgebra generated by all the lower degree classes (and up to a coefficient that is of no importance in characteristic $0$) the
classes $s_i$ and $c_i$ are often interchangeable when studying the differentials in the spectral sequence.

\end{remark}

\begin{nparagraph}
From this point forward we restrict the discussion to $\V_{n,1}$. From the previous lemma we see that the space of invariants
in $\Lambda^\bullet(V) \tensor \Lambda^\bullet(S^2(V^*) \tensor V)$ is generated by the classes $c_i$ of degree $2i$, for $1 \le i \le n$,
with no relations imposed on them.

Let $BGL(n)$ be the classifying space of the topological group $GL(n, \CC)$ and
$$
p\from EGL(n) \to BGL(n)
$$
be the tautological principal
$GL(n)$-bundle. Recall that the de Rham cohomology
$$
H^\bullet(BGL(n), \k) \isom \k[c_1, c_2, \ldots, c_n],
$$
with generators $c_i$ of degree $2i$.
Consider $2d$-dimensional skeleton $\sk_{2d}BGL(n) \into BGL(n)$ and denote $X_{2d} = p^{-1}(\sk_{2d}BGL(n))$. We will also denote
by $\SS$ the topological suspension functor. Now we are ready to prove the main theorem of this paper.
\end{nparagraph}

\begin{theorem}
We have an isomorphism
$$
H^\bullet(\V_{n,1}, \k) \ \isom\  H^\bullet(\SS^2 X_{2(n-1)}, \k).
$$
\end{theorem}

\proof
Consider the spectral sequence associated to the Lie subalgebra $\gl(n,1) \into \V_{n,1}$. From the discussion above and lemma \ref{lemma_invar}
we find that $E_1^{pq} = 0$ if $p$ is odd and
$$
E_1^{2p,q} = \bigoplus_{|\lambda| = p} H^\bullet(\gl(n,1), \Sigma^{\wtilde\lambda}(V) \tensor \Sigma^{\wtilde\lambda}(V^*)).
$$
Next, from proposition \ref{prop_cohom_gln1} we see that
$$
E_1^{2p,q} = \begin{cases}
\bigoplus_{|\lambda| = p} \k[e_1, \ldots, e_{2n-1}],&\text{if $p \le n - 1$,}\\
\bigoplus_{|\lambda| = p} \k[e_1, \ldots, e_{2n-1}, e'_1],&\text{if $p \ge n$.}
\end{cases}
$$

We will also consider the spectral sequence for the cohomology of the total space of the tautological bundle $EGL(n) \to BGL(n)$ with coefficients
in $\k[e'_1]$
$$
F_1^{\bullet\bullet} = \k[e_1, \ldots, e_{2n-1}] \tensor \k[c_1, \ldots, c_n] \tensor \k[e'_1] \Rightarrow H^\bullet(EGL(n), \k[e'_1]) \isom \k[e'_1],
$$
since the $EGL(n)$ is contractible.

Let us look at the short exact sequence $E_1 \into F_1 \epi E_1 / F_1$. For the quotient we have
$$
(E_1 / F_1)^{\bullet\bullet} = \left(\k[e_1, \ldots, e_{2n-1}] \tensor \k[c_1, \ldots, c_n] / I_n \right) e'_1,
$$
where $I_n$ is spanned by all the monomials of total degree at least $2n$. The standard argument shows that it converges to
$H^\bullet(X_{2(n-1)}, \k)\cdot e'_1$. Now from the long exact sequence of cohomology associated with this triple of spectral sequences
we find that
\begin{align*}
H^0(\V_{n,1}, \k) &= \k,\\
H^1(\V_{n,1}, \k) &= H^2(\V_{n,1}, \k) = 0,\\
H^i(\V_{n,1}, \k) &\isom H^{i-2}(X_{2(n-1)},\k), \quad \text{for $i > 2$}.
\end{align*}

\qed

\begin{remark}[Cohomology of $\V_{1,n}$.]
As another application of the methodology developed here we give an alternative proof of
$$
H^\bullet(\V_{1,n}, \k) \ \isom\ H^\bullet(S^{2n-1}, \k),
$$
which is of course a special case of theorem \ref{athm_mlessn}.

As before we write $V$ for the standard representation of $\gl(n,1)$ and let $W$ be the standard representation of $\gl(1,n)$. Notice
that they differ by the change of parity, i.e. $W = \Pi V$. Therefore, up to parity we have isomorphisms $\Sigma^\lambda(W) \isom \Sigma^{\lambda'}(V)$.
Now, from lemma \ref{lemma_invar} and the discussion in remark \ref{rem_invar} we find that diagrams $\lambda$ contributing to the first layer of the
spectral sequence are contained in the thick hook with $1$ row and $n$ columns. Consider a free algebra generated by elements $s_i$ of degree $2i$,
for all $i \ge 1$ and the ideal $J = (s_j \mid j \ge n)$. Then we see that the space of invariants can be described as
$$
\left(\Lambda^\bullet(W) \tensor \Lambda^\bullet(S^2(W^*) \tensor W)\right)^{\gl(1,n)} \ \isom\ \k[s_1, s_2, \ldots ] / (J)^2.
$$

Denote by $z_\lambda$ the basis vector in the first layer of the spectral sequence coming from $\Sigma^\lambda(W) \tensor \Sigma^\lambda(W^*)$.
From proposition \ref{prop_cohom_gln1} and the relation between $W$ and $V$ mentioned above we see that elements $e'_1$ appear only for
$z_\lambda$ with $\rmht(\lambda') \ge n$. Proceeding as in the proof of proposition \ref{prop_cohom_gln1} we consider decreasing filtration
$\cF_l E_2$ by $l_n(\lambda')$ (see \ref{prop_proof_n-1} for the notation) and another grading on $E_2$ by the number of $e'_1$, which makes
it a collection of bicomplexes.
After passing to the associated graded one of the differentials in the bicomplex can be written as $d(e'_1 z_\lambda) = z_{\lambda_+}$, where
$\lambda_+ = ((\lambda')^+)'$, i.e. obtained from $\lambda$ by adding one box in the first row.
Taking cohomology with respect to this differential (and keeping the other differential intact) we obtain a spectral sequence
which starting from the second layer is isomorphic to another spectral sequence $F_\bullet$,
that starts with the first layer isomorphic to
$$
F_1 \isom \k[e_1, \ldots, e_{2n-1}] \tensor \k[s_1, \ldots, s_{n-1}].
$$
The elements $e_{2i-1}$ for $1 \le i < n$ are transgressive, i.e. they are killed by differentials of the spectral sequence, hence we find
$$
H^\bullet(\V_{1,n}, \k) \ \isom\ \k[e_{2n-1}].
$$

\qed

\end{remark}

\vfill\eject

\end{document}